\documentclass[12pt]{article}
\usepackage{amssymb}
\usepackage{arydshln} 
\newtheorem{const}{Construction}
\newtheorem{prop}{Proposition}

\begin{document}

\begin{center}
\Large
The Heritage of Cayley-Sudoku Tables\footnote{Portions of this paper appear in  Kady Hossner Boden \& Michael B. Ward (2019) A New Class of Cayley-Sudoku Tables, Mathematics Magazine, 92:4, 243-251, DOI: 10.1080/0025570X.2019.1613949 }
\end{center}

\begin{flushright}
Kady Hossner Boden\\
St. Stephen's Academy \\
Beaverton, OR 97008 \\
\verb+kboden@ststephensacademy.com+ \\

\vskip 2 mm

Michael B. Ward\\
Western Oregon University \\
Monmouth, OR 97361\\
\verb+wardm@wou.edu+
\end{flushright}

\begin{abstract}
A Cayley-Sudoku table of a finite group G is a Cayley  table for G subdivided into uniformly sized rectangular blocks, in such a way that each group element appears once in each block.  Obviously inspired by the popularity of Sudoku puzzles, Cayley-Sudoku tables and three ways to construct them were introduced by J. Carmichael, K. Schloeman, and M. B. Ward [\ref{csw}].  Since then, we have discovered that the first two constructions have an unexpected heritage in the work of two distinguished mathematicians, Reinhold Baer and Jozsef Denes.  Special cases of the first construction are reinvented in recent publications.  This paper has four aims.  First, we review Constructions 1 and 2 and uncover their heritage.  Next we turn to some new instances of Construction 2 inspired by Baer, which answer an open question in [\ref{csw}].  Third, we provide a very brief outline of the recent appearances of Construction 1.   We conclude with an invitation to seek out the heritage of Construction 3.
\end{abstract}

\paragraph*{Introduction} A Cayley-Sudoku table of a finite group $G$ is a Cayley (i.e. operation) table for $G$ subdivided into uniformly sized rectangular blocks, in such a way that each group element appears once in each block. For example, Table \ref{z9} is a Cayley-Sudoku table for $\mathbb{Z}_9 := \{0,1,2,3,4,5,6,7,8\}$ under addition modulo 9 and Table \ref{s3} is a Cayley-Sudoku table for $S_3$, the symmetric group on three symbols.  Obviously inspired by the popularity of Sudoku puzzles, Cayley-Sudoku tables were introduced in [\ref{csw}], which gave three constructions.  The second construction involved a curious condition.  After several fruitless inquiries, M. Ward asked about that condition in a talk at the XXX Ohio State-Denison Math Conference where, fortuitously, quasigroup theorists were in attendance.  Numerous hands went up.  First responder Professor Clifton E. Ealy, Jr.\footnote{Western Michigan University} announced ``You and your students have rediscovered a 1939 theorem of Reinhold Baer!'' Having been pointed in the right direction, we now know Constructions 1 and 2 of the three constructions in [\ref{csw}] have an unexpected heritage in the work of two distinguished mathematicians, Baer [\ref{baer}] and J\'{o}zsef D\'{e}nes [\ref{denes1}, \ref{denes2}].  Camouflaged special cases of Construction 1 also appear in more recent publications [\ref{pv},\ref{lorch}], adding to its heritage.

This paper has four aims.  First, we review Constructions 1 and 2 and uncover their heritage.  The approach is descriptive and rather informal, assuming only a familiarity with elementary group theory through cosets.  Next we turn to some new instances of Construction 2 inspired by Baer, which are more technical and require some experience with permutation groups.  These instances answer an open question in [\ref{csw}].  Third, we provide a very brief outline of recent appearances of Construction 1.   We conclude with an invitation to seek out the heritage of Construction 3.

\begin{table}[h!]
\begin{center}
\begin{tabular}{c||ccc|ccc|ccc|}
\null & 0 & 3 & 6 & 1 & 4 & 7 & 2 & 5 & 8 \\
\hline\hline 0 & 0 & 3 & 6 & 1 & 4 & 7 & 2 & 5 & 8 \\
1 & 1 & 4 & 7 & 2 & 5 & 8 & 3 & 6 & 0 \\
2 & 2 & 5 & 8 & 3 & 6 & 0 & 4 & 7 & 1 \\
\hline 3 & 3 & 6 & 0 & 4 & 7 & 1 & 5 & 8 & 2 \\
4 & 4 & 7 & 1 & 5 & 8 & 2 & 6 & 0 & 3 \\
5 & 5 & 8 & 2 & 6 & 0 & 3 & 7 & 1 & 4 \\
\hline 6 & 6 & 0 & 3 & 7 & 1 & 4 & 8 & 2 & 5 \\
7 & 7 & 1 & 4 & 8 & 2 & 5 & 0 & 3 & 6 \\
8 & 8 & 2 & 5 & 0 & 3 & 6 & 1 & 4 & 7 \\ \hline
\end{tabular}
\end{center}
\caption{$\mathbb{Z}_9$ Cayley-Sudoku Table} \label{z9}
\end{table}

\begin{table}[h!]
\begin{center}
\begin{tabular}{c||ccc|ccc|}
\null & $(1)$ & $(13)$ & $(132)$ & $(12)$ & $(123)$ & $(23)$   \\ \hline \hline
$(1)$ & $(1)$ & $(13)$ & $(132)$ & $(12)$ & $(123)$ & $(23)$ \\
$(12)$ & $(12)$ & $(123)$ & $(23)$ & $(1)$ & $(13)$ & $(132)$ \\ \hline
$(13)$ & $(13)$ & $(1)$ & $(12)$ & $(132)$ & $(23)$ & $(123)$ \\
$(132)$ & $(132)$ & $(23)$ & $(123)$ & $(13)$ & $(1)$ & $(12)$ \\ \hline
$(23)$ & $(23)$ & $(132)$ & $(13)$ & $(123)$ & $(12)$ & $(1)$ \\
$(123)$ & $(123)$ & $(12)$ & $(1)$ & $(23)$ & $(132)$ & $(13)$ \\ \hline
\end{tabular}
\end{center}
\caption{$S_3$ Cayley-Sudoku Table} \label{s3}
\end{table}

\paragraph*{Constructing Cayley-Sudoku Tables}
Tables \ref{z9} and \ref{s3} illustrate Construction 1 of [\ref{csw}].  Consider the subgroup $S = \left< 3 \right> = \{0,3,6\}$ of $\mathbb{Z}_9$. Columns in each block are labeled by elements of the right cosets $S+0 = \{0,3,6\}$, $S+1=\{1,4,7\}$, and $S+2=\{2,5,8\}$.  The rows in each block are labeled with a complete set of left coset representatives of $S$ (that is, one element from each left coset), $L_1 = \{0,1,2\}$, $L_2 = \{3,4,5\}$, and $L_3 = \{6,7,8\}$.  (The distinction between right and left cosets is unimportant in the commutative group $\mathbb{Z}_9$, but it is critical in general.)

In the sequel, we usually economize by indicating the row and column labels for the blocks by listing the sets of labels rather than individual labels.  Thus, the layout of Table \ref{z9} is

\begin{center}
\begin{tabular}{c||c|c|c|}
\null & $S+0$ & $S+1$ & $S+2$ \\ \hline \hline
$L_1$ & & & \\ \hline
$L_2$ & & & \\ \hline
$L_3$ & & & \\ \hline
\end{tabular}
\end{center}

For $S_3$, consider the subgroup $S=\left< (12) \right>$.  Now the rows in each block are labeled with elements of the left cosets $(1)S=\{(1),(12)\}$, $(13)S=\{(13),(132)\}$, and $(23)S=\{(23),(123)\}$.\footnote{Permutations here are composed left to right.} The right cosets are $S(1)=\{(1),(12)\}$, $S(13)=\{(13),(123)\}$, and $S(23)=\{(23),(132)\}$ and the columns in each block are labeled with a complete set of right coset representatives of $S$, $R_1=\{(1),(13),(132)\}$ and $R_2=\{(12),(123),(23)\}$.  The condensed layout of Table \ref{s3} is

\begin{center}
\begin{tabular}{c||c|c|}
\null & $R_1$ & $R_2$  \\ \hline \hline
$(1)S$ & & \\ \hline
$(13)S$ & & \\ \hline
$(23)S$ & & \\ \hline
\end{tabular}
\end{center}

Construction 1 simply says either such layout always produces a Cayley-Sudoku table.

\begin{const} \label{const1} Let $G$ be a finite group. Assume $S$ is a subgroup of $G$ having order $k$ and index $n$.  If $Sg_{1}, Sg_{2},\ldots,Sg_{n}$ are the distinct right cosets of $S$ in $G$, then arranging the Cayley table of $G$ with columns labeled by the cosets $Sg_{1}, Sg_{2},\ldots,Sg_{n}$ and the rows labeled by sets $L_{1}, L_{2}, \ldots, L_{k}$ (as in Table \ref{const1}) yields a Cayley-Sudoku table of $G$ with blocks of dimension $n\times k$ if and only if $L_{1}, L_{2}, \ldots, L_{k}$ partition G into complete sets of left coset representatives of $S$ in $G$.

\begin{table}[h]
\begin{center}\begin{tabular}{c||c|c|c|c|}
  & $Sg_{1}$ & $Sg_{2}$ & $\ldots$ & $Sg_{n}$ \\
\hline \hline $L_{1}$ &  &  &  &  \\
\hline $L_{2}$ &  &  &  &  \\
\hline $\vdots$ &  &  &  &  \\
\hline $L_{k}$ &  &  &  &  \\
\hline
\end{tabular} \end{center}
\caption{Construction \ref{const1}R Using Right Cosets and Left Coset Representatives} \label{const1R}
\end{table}

Furthermore, if $y_{1}S, y_{2}S,\ldots,y_{n}S$ are the $n$ distinct left cosets of $S$ in $G$, then arranging the Cayley table of $G$ with rows labeled by the cosets $y_{1}S, y_{2}S,\ldots,y_{n}S$ and the columns labeled by sets $R_{1}, R_{2}, \ldots, R_{k}$ yields a Cayley-Sudoku table of $G$ with blocks of dimension $k \times n$ if and only if $R_{1}, R_{2}, \ldots, R_{k}$ partition G into complete sets of right coset representatives of $S$ in $G$.

\begin{table}[h]
\begin{center}\begin{tabular}{c||c|c|c|c|}
  & $R_1$ & $R_2$ & $\ldots$ & $R_k$ \\
\hline \hline $y_{1}S$ &  &  &  &  \\
\hline $y_{2}S$ &  &  &  &  \\
\hline $\vdots$ &  &  &  &  \\
\hline $y_{n}S$ &  &  &  &  \\
\hline
\end{tabular} \end{center}
\caption{Construction \ref{const1}L Using Left Cosets and Right Coset Representatives} \label{const1L}
\end{table}

\end{const}

Notice the second part is dual to the first, obtained by switching right with left and rows with columns.  We refer to the first as Construction 1R because it uses right cosets and the second one as 1L.  For any subgroup of any finite group, one can always partition the group into complete sets of left or right coset representatives.  Therefore, every group has Cayley-Sudoku tables corresponding to each of its subgroups.

We now turn to the next construction, for which we reveiw a definition.  If $S$ is a subgroup of the group $G$ and $g \in G$, $S^g$ denotes the subgroup $g^{-1}Sg:=\{g^{-1}sg : s \in S\}$, which is called a \emph{conjugate} of $S$.  

\begin{const} \label{const2}  Assume $S$ is a subgroup of $G$ having order $k$ and index $n$.  Also suppose $y_{1}S$, $y_{2}S$, \ldots, $y_{n}S$ are the distinct left cosets of $S$ in $G$. Arranging the Cayley table of $G$ with columns labeled by the cosets $y_{1}S$, $y_{2}S$, \ldots, $y_{n}S$ and the rows labeled by sets $L_{1}, L_{2}, \ldots, L_{k}$ yields a Cayley-Sudoku table of $G$ with blocks of dimension $n \times k$ if and only if $L_{1}, L_{2}, \ldots, L_{k}$ are complete sets of left coset representatives of $S^{g}$ for all $g \in G$. \\

\begin{table}[h]
\begin{center} \begin{tabular}{c||c|c|c|c|}
 & $y_{1}S$ & $y_{2}S$ & $\ldots$ & $y_{n}S$   \\
\hline \hline $L_{1}$ &  &  &  &    \\
\hline $L_{2}$ &  &  &  &    \\
\hline $\vdots$ &  &  &  &  \\
\hline $L_{k}$ &  &  &  &  \\ \hline
\end{tabular} \end{center}
\caption{Construction \ref{const2}L Using Left Cosets and Left Coset Representatives} \label{const2L}
\end{table}

Furthermore, suppose $Sg_1$, $Sg_2$, \ldots, $Sg_n$ are the distinct right cosets of $S$ in $G$.  Arranging the Cayley table of $G$ with rows labeled by the cosets $Sg_1$, $Sg_2$, \ldots, $Sg_n$ and columns labeled by the sets $R_1$, $R_2$, \ldots $R_k$ yields a Cayley-Sudoku table of $G$ if and only if $R_1$, $R_2$, \ldots $R_k$ partition $G$ into complete sets of right coset representatives of $S^g$ for all $g \in G$.
\begin{table}[h!]
\begin{center}
\begin{tabular}{c||c|c|c|c|}
\null & $R_0$ & $R_1$ & $\cdots$ & $R_k$ \\ \hline \hline
$Sg_1$ & & & & \\ \hline
$Sg_2$ & & & & \\ \hline
\vdots & & & & \\ \hline
$Sg_m$ & & & & \\ \hline
\end{tabular}
\end{center}
\caption{Construction \ref{const2}R Using Right Cosets and Right Coset Representatives} \label{const2R}
\end{table}
\end{const}

We refer to the two parts as Construction 2L and 2R, respectively, referring again to the use of left or right cosets.

To illustrate Construction 2L, let $G = S_3$ and $S= \left< (1,2) \right>$.  The left cosets are $(1)S=\{(1),(12)\}$, $(13)S=\{(13),(132)\}$, and $(23)S=\{(23),(123)\}$.  The conjugates of $S$ in $G$ are $\left< (12) \right>$, $\left< (13) \right>$, and $\left< (23) \right>$.  It is easy to check that $L_1 = \{(1), (123), (132) \}$ and $L_2 = \{(12), (13), (23) \}$ partition $G$ into complete sets of left coset representatives for each of the conjugates of $S$ in $G$.  Thus, table \ref{s32} is an instance of Construction 2L.

\begin{table}[h!]
\begin{center}
\begin{tabular}{c||cc|cc|cc|}
\null & $(1)$ & $(12)$ & $(13)$ & $(132)$ & $(23)$ & $(123)$   \\ \hline \hline
$(1)$ & $(1)$ & $(12)$ & $(13)$ & $(132)$ & $(23)$ & $(123)$ \\
$(123)$ & $(123)$ & $(23)$ & $(12)$ & $(1)$ & $(13)$ & $(132)$ \\
$(132)$ & $(132)$ & $(13)$ & $(23)$ & $(123)$ & $(12)$ & $(1)$ \\ \hline
$(12)$ & $(12)$ & $(1)$ & $(123)$ & $(23)$ & $(132)$ & $(13)$ \\
$(13)$ & $(13)$ & $(132)$ & $(1)$ & $(12)$ & $(123)$ & $(23)$ \\
$(23)$ & $(23)$ & $(123)$ & $(132)$ & $(13)$ & $(1)$ & $(12)$ \\
\end{tabular}
\end{center}
\caption{Another $S_3$ Cayley-Sudoku Table} \label{s32}
\end{table}

Construction 2L looks very similar to Construction 1R, but 1R required \emph{right} cosets and left coset representatives.  The price we pay for using \emph{left} cosets along with left coset representatives in 2L is that the coset representatives must be complete sets of left coset representatives not just for the subgroup $S$ but for all the conjugates of $S$ at once.  This can be a high price.  For example, it is impossible to find such representatives for the subgroup $\left<(12)(34)\right>$ in the symmetric group $S_4$.  The obvious question ``When \emph{can} we get such coset representatives?" led to the enthusiastic reference to Baer's Theorem mentioned in the introduction.

The final construction shows how to extend a Cayley-Sudoku table of a subgroup to a Cayley-Sudoku table of the full group.  See [\ref{csw}] for an example of its use.

\begin{const} \label{const3} Let $G$ be a finite group with a subgroup $A$.  Let $C_{1},C_{2}, \ldots, C_{k}$ partition $A$ and $R_{1},R_{2}, \ldots R_{n}$ partition $A$ such that the following table is a Cayley-Sudoku table of $A$.
\begin{center}
\begin{tabular}{c||c|c|c|c|}
 & $C_{1}$ & $C_{2}$ & $\ldots$ & $C_{k}$ \\
\hline $R_{1}$ &  &  &  &  \\
\hline $R_{2}$ &  &  &  &  \\
\hline $\vdots$ &  &  &  &  \\
\hline $R_{n}$ &  &  &  &  \\
\hline
\end{tabular} \end{center}
If $\lbrace l_{1},l_{2}, \ldots, l_{t} \rbrace$ and $\lbrace r_{1},r_{2},\ldots r_{t}\rbrace$ are complete sets of left and right coset representatives, respectively, of $A$ in $G,$ then arranging the Cayley table of $G$ with columns labeled with the sets $C_{i}r_{j}$, $i=1, \ldots, k$, $j=1, \ldots , t$ and the $b^{th}$ block of rows labeled with $l_jR_b$, $j=1, \ldots , t$, for $b = 1, \ldots, n$ (as in Table \ref{const3Table}) yields a Cayley-Sudoku table of $G$ with blocks of dimension $tk \times n$.
\end{const}

\begin{table}[h!]
\begin{center}
\begin{tabular}{c||c|c|c|c|c|c|c|c|c|c|c|}
& $C_{1}r_{1}$ & $C_{2}r_{1}$ & $\ldots$ & $C_{k}r_{1}$ & $C_{1}r_{2}$ & $\ldots$ & $C_{k}r_{2}$ & $\ldots$ & $C_{1}r_{t}$ & $\dots$ & $C_{k}r_{t}$ \\ \hline \hline
$l_{1}R_{1}$ &  &  &  &  &  &  &  &  &  &  &  \\
$l_{2}R_{1}$ &  &  &  &  &  &  &  &  &  &  &  \\
$\vdots$ &  &  &  &  &  &  &  &  &  &  &  \\
$l_{t}R_{1}$ &  &  &  &  &  &  &  &  &  &  &  \\ \hline
$l_{1}R_{2}$ &  &  &  &  &  &  &  &  &  &  &  \\
$\vdots$ &  &  &  &  &  &  &  &  &  &  &  \\
$l_{t}R_{2}$ &  &  &  &  &  &  &  &  &  &  &  \\ \hline
$\vdots$ &  &  &  &  &  &  &  &  &  &  &  \\  \hline
$l_{1}R_{n}$ &  &  &  &  &  &  &  &  &  &  &  \\
$\vdots$ &  &  &  &  &  &  &  &  &  &  &  \\
$l_{t}R_{n}$ &  &  &  &  &  &  &  &  &  &  &  \\ \hline
\end{tabular} \end{center}
\caption{Construction \ref{const3}} \label{const3Table}
\end{table}

\paragraph*{The Heritage of Construction 1}
Construction 1 is a rediscovery and clarification of a theorem of D\'{e}nes [\ref{denes1}, \ref{denes2}].  An \emph{$(m,1)$-complete Latin rectangle} is a rectangle that can be completed to a Latin square and contains $m$ different symbols each occurring exactly once. Since every Cayley table is a (bordered) Latin square, the blocks in our Cayley-Sudoku tables are $(m,1)$-complete Latin rectangles where $m$ is the order of the group.  D\'{e}nes stated

\textsc{Theorem} [\ref{denes2}, Theorem 1.5.5] \emph{If $L$ is the Latin square representing the multiplication [Cayley] table of a group $G$ of order $m$, where $m$ is a composite number, then $L$ can be split [partitioned] into a set of $m$ $(m,1)$-complete non-trivial\footnote{Not consisting of a single row or column} Latin rectangles.}

In proving the theorem, D\'{e}nes takes a proper non-trivial subgroup (which exists since $m$ is composite) of $G$ and arranges the Cayley table exactly as specified in Construction 1L--with one possible flaw.  The coset representatives he designates as column labels might not be right coset representatives as required.  Use of right versus left cosets is ambiguous in the proof and the examples in [\ref{denes1}] and [\ref{denes2}] use normal subgroups where the distinction is irrelevant.  Nevertheless, D\'{e}nes gets credit for constructing the first Cayley-Sudoku table.  The authors of [\ref{csw}] can console themselves, however, for republishing a known result with the thought that their proof might be clearer.

(By the way, the use of a proper non-trivial subgroup is only to ensure the non-triviality of the Latin rectangles or blocks, as we call them. Since we allow blocks consisting of a single rows or columns, we omit the hypothesis that the order of the group is composite.)

\paragraph*{The Heritage of Construction 2}

A set $Q$ with a binary operation $\cdot$ is a  \emph{quasigroup} provided its Cayley table is a bordered Latin square on the elements of $Q$ or, equivalently, it has cancelation, which is to say, for every $x, y, a \in Q$, if $a \cdot x = a \cdot y$, then $x=y$ (row $a$ in the Cayley table contains each element only once) and $x \cdot a = y \cdot a$ implies $x=y$ (column $a$ in the Cayley table contains each element only once).  Removing the borders of Table \ref{loop}, for example, leaves a Latin square, so it gives a quasigroup.  In a quasigroup, there need not be inverses, an identity, or even associativity.

With $S$ and $G$ as in Construction 2R, fix a complete set of right coset representatives $\{r_1, r_2, \ldots, r_m\}$ of $S$ in $G$.  (Think of this as one of the $R_i$ in Construction 2.)  Define a coset multiplication by $Sr_i \cdot Sr_j := Sr_i r_j$.  Because the coset representatives are fixed, this gives a binary operation on $\mathcal{R}$, the set of right cosets of $S$ in $G$.  Normally, coset multiplication like that is well-defined only when $S$ is a normal subgroup of $G$.  However, by fixing the coset representatives in advance, all is well.  We can now state the 1939 theorem mentioned in the introduction.
\vspace{1ex}

\textbf{Baer's Theorem} \label{baerthm} [\ref{baer}, Theorem 2.3] \emph{$\mathcal{R}$ under $\cdot$ as defined above is a quasigroup\footnote{Baer calls it a division system.} if and only if $\{r_1, r_2, \ldots, r_m\}$ is a complete set of right coset representatives of $S^g$ for every $g \in G$.}
\vspace{1ex}

That last bit is the same condition as in Construction 2R.  Therefore, we have
\vspace{1ex}

\textbf{Construction 2R \`{a} la Baer} \emph{With notation as in Construction \ref{const2}R, the following are equivalent.}

\emph{(a) The arrangement of the Cayley table in Table \ref{const2R} gives a Cayley-Sudoku table.
}

\emph{(b) The sets $R_1$ through $R_k$ partition $G$ and each is a complete set of right coset representatives of $S^g$ for every $g \in G$.
}

\emph{(c) The sets $R_1$ through $R_k$ partition $G$ and each gives rise to a quasigroup on the right cosets of $S$ as described above.}
\vspace{1ex}

We chose Construction \ref{const2}R in order to conform with Baer's use of right cosets.  One can also prove the left-handed version corresponding to Construction 2L.

It may be of interest to prove this new formulation of Construction \ref{const2} without an appeal to Baer's Theorem, thereby providing a new proof of Baer's Theorem.  It suffices to prove (a) is equivalent to (c).

 Consider a table $T$ laid out as in Table \ref{const2R}.  Take any $R_j = \{r_1, r_2, \ldots, r_m\}$ and consider the column of blocks underneath it.  Relabeling the elements of $R_j$ if necessary, we may assume $Sg_i = Sr_i$ for each $i$.  Substituting, writing out the individual column elements in $R_j$, and filling-in the blocks, the column of blocks headed by $R_j$ is transformed as follows.
\begin{center}
\begin{tabular}{c||c|}
\null & $R_j$ \\ \hline \hline
$Sg_1$ & \\ \hline
$Sg_2$ & \\ \hline
\vdots & \\ \hline
$Sg_m$ & \\ \hline
\end{tabular}
$\longrightarrow$
\begin{tabular}{c||cccc|}
\null & $r_1$ & $r_2$ & $\cdots$ & $r_m$ \\ \hline \hline
$Sr_1$ & $Sr_1 r_1$ & $Sr_1 r_2$ & $\cdots$ & $Sr_1 r_m$ \\ \hline
$Sr_2$ & $Sr_2 r_1$ & $Sr_2 r_2$ & $\cdots$ & $Sr_2 r_m$ \\ \hline
\vdots & & & &\\ \hline
$Sr_m$ & $Sr_m r_1$ & $Sr_m r_2$ & $\cdots$ & $Sr_m r_m$ \\ \hline
\end{tabular}
\end{center}

Let $L$ denote the transformed column of blocks without its borders. At this point, under our convention, each coset appearing in $L$ is really just used as shorthand for the (vertical) list of its elements. However, we can also think of the cosets simply as symbols in the square array $L$.  Exploiting those two points of view, we claim for each $j$, $L$ is a Latin square on the $m$ right cosets of $S$ if and only if the original table $T$ is a Cayley-Sudoku table.  It is easy to see that the cosets listed in any column of $L$ are distinct since $Sr_1, \ldots, Sr_m$ are distinct.  Turning to the rows, recall each row in the transformed table is actually one block of $T$ (thinking now of the cosets as lists).  Thus, the sudoku condition that each group element appear exactly once in each block of $T$ is equivalent to the cosets listed in each row of $L$ being distinct, establishing our claim.

Reintroduce the borders on $L$, but change the top border from $r_1, r_2, \ldots, r_m$ to $Sr_1$, $Sr_2$, \ldots, $Sr_m$ (and regard the cosets as symbols).  The resulting table (see below) is precisely the table for Baer's coset multiplication based on the right coset representatives in $R_j$ and, by the previous paragraph, it defines a quasigroup if and only if $T$ is a Cayley-Sudoku table, as was to be shown.

\begin{center}
\begin{tabular}{c||cccc|}
\null & $Sr_1$ & $Sr_2$ & $\cdots$ & $Sr_m$ \\ \hline \hline
$Sr_1$ & $Sr_1 r_1$ & $Sr_1 r_2$ & $\cdots$ & $Sr_1 r_m$ \\ \hline
$Sr_2$ & $Sr_2 r_1$ & $Sr_2 r_2$ & $\cdots$ & $Sr_2 r_m$ \\ \hline
\vdots & & & &\\ \hline
$Sr_m$ & $Sr_m r_1$ & $Sr_m r_2$ & $\cdots$ & $Sr_m r_m$ \\ \hline
\end{tabular}
\end{center}

\paragraph*{New Instances of Construction 2}

 Other than the trivial case where $S$ is a normal subgroup (and Construction 2 reduces to Construction 1), the following two theorems gave the only general setting known to the authors of [\ref{csw}] wherein Construction 2 could be applied.  The authors asked for other such settings.

\begin{prop} \label{need1} Assume $S$ is a subgroup of a finite group $G$.  Suppose $R$ is a complete set of right [left] coset representatives of $S^g$ for all $g \in G$.  Then the sets $sR := \{sr : r \in R\}$ [$Rs := \{rs : r \in R\}$],  $s \in S$ partition $G$ into complete sets of right [left] coset representatives of $S^g$ for all $g \in G$.
\end{prop}

In other words, in applying Construction \ref{const2}, it is sufficient to find one set of coset representatives of the desired sort.

\begin{prop} \label{complement} Suppose $S$ is a subgroup of the finite group $G$ and there is a subgroup $C$ such that $G=CS$ and $C \cap S = 1$ (i.e. $C$ is a \emph{complement} for $S$), then $C$ is a complete set of left and right coset representatives of $S$ in $G$.
\end{prop}

Table \ref{s32} is an instance of Propositions \ref{complement} and \ref{need1} since $L_1$ is a complement for $S$ in $G$ and $L_2 = L_1(12)$.

While Baer's Theorem gives another way to think about Construction 2, it does not directly give the new instances of the construction called for in [\ref{csw}], ``new" meaning not accounted for by Proposition \ref{complement}.  Nevertheless, reading Baer inspired a new class of examples arising from quasigroups, which we we will now explain, beginning with a quick overview of some basics about permutation groups and quasigroups.

Suppose $G$ is a group of permutations of a set $A$.
For each $a \in A$ and $g \in G$, $a^g$ denotes the image of $a$ under $g$ and $G_a := \{g \in G : a^g = a \}$ denotes the \emph{stabilizer} of $a$ in $G$.\footnote{We beg the reader's pardon for  writing $a^g$ in place of the more familiar function notation $g(a)$.  It is common practice among group theorists,  and  it is necessary in order to conform to upcoming  GAP [\ref{gap}] calculations.}
$G$ is \emph{transitive} when for every $a, b \in A$ there is a $g \in G$ such that $a^g = b$.
$G$ is \emph{regular} provided $G$ is transitive and $G_a = 1$ for any $a \in A$.
The following is a standard result.

\begin{prop} \label{permlemma}
Suppose $G$ is a transitive group of permutations on a set $A$ and $a \in A$.

(a) For each $g \in G$, $(G_a)^g := g^{-1} G_a g = G_{a^g}$.

(b) If $T$ is a regular subgroup of $G$, then $|T| = |A|$.

(c) If $T$ is a complement of $G_a$ in $G$, then $T$ is regular.
\end{prop}

For the convenience of the reader, we prove the next ``well-known'' proposition.

\begin{prop} \label{an} Suppose $n$ is a positive integer and $n \equiv 2 (\textrm{mod } 4)$, then the alternating group $A_n$ does not contain a regular subgroup.
\end{prop}

Proof:  With $n$ as in the hypotheses, assume $T$ is a regular subgroup of $A_n$.  By Proposition \ref{permlemma}, $|T| = n$.  Thus, $T$ contains an element $t$ of order 2 by Cauchy's Theorem.  Since $T$ is regular, $t$ has no fixed points.  Therefore, it is the composition of $n/2$ 2-cycles, a contradiction since $n/2$ is odd.
\vspace{1ex}


Suppose $Q$ with operation $\cdot$ is a finite quasigroup.
For each $\ell \in Q$ define $\rho_\ell, \lambda_\ell: Q \rightarrow Q$  by $q^{\rho_\ell} = q \cdot \ell$ and $q^{\lambda_\ell} = \ell \cdot q$. 

Using the definition of quasigroup, one can easily prove the next result.

\begin{prop} \label{quasilemma} Suppose $Q$ is a finite quasigroup.  For each $\ell \in Q$,  $\rho_\ell$ and $\lambda_\ell$ are permutations of $Q$.  For every $a, b \in Q$ there exist $\ell, l \in Q$ such that $a^{\rho_\ell} = b$ and $a^{\lambda_l} = b$.
\end{prop}

When $Q$ is a group, $\{\rho_\ell : \ell \in Q \}$ and $\{\lambda_\ell : \ell \in Q \}$ are permutation groups isomorphic to $Q$.  That is Cayley's Representation Theorem.  In general, however, those sets are not groups.  Instead, quasigroup theorist consider the groups of permutations generated by those sets, $RMult(Q) := \left< \rho_\ell: \ell \in Q  \right>$ and $LMult(Q) := \left< \lambda_\ell: \ell \in Q  \right>$.  By Proposition \ref{quasilemma}, each of these groups is transitive.

We can now describe instances of Construction \ref{const2} arising from quasigroups.

\begin{prop} \label{quasieg} Suppose $Q$ is a quasigroup and $c \in Q$. Let $G = RMult(Q)$, then $R := \{\rho_\ell : \ell \in Q \}$ is a complete set of right coset representatives of $(G_c)^g$ for every $g \in G$.  Similarly, if $G = LMult(Q)$, then and $L := \{\lambda_\ell : \ell \in Q \}$ is a complete set of left coset representatives of $(G_c)^g$ for every $g \in G$.
\end{prop}

Proof.  Let $G = RMult(Q)$ and $g, h \in G$.  By Proposition \ref{quasilemma}, there exists $\ell \in Q$ such that $(c^g)^{\rho_\ell} = c^{gh}$.  Therefore, $\rho_\ell h^{-1} \in G_{c^g}$ which equals $(G_c)^g$ by Proposition \ref{permlemma}.  Thus, $(G_c)^g \rho_\ell = (G_c)^g h$ and $R$ contains a representative of each coset of $(G_c)^g$.  Moreover, for $\ell, m \in Q$, $(G_c)^g \rho_\ell = (G_c)^g \rho_m$ if and only if $(c^g)^{\rho_\ell} = (c^g)^{\rho_m}$, that is, $(c^g) \cdot \ell = (c^g) \cdot m$.  Therefore, $\ell = m$ by cancelation in $Q$.  Therefore, $R$ is a complete set of right coset representatives of $(G_c)^g$ for every $g \in G$.  The proof for $G = LMult(Q)$ is similar.
\vspace{1ex}

Proposition \ref{quasieg} and Proposition \ref{need1} imply that for any quasigroup $Q$, Construction \ref{const2} applies to $G = LMult(Q)$ and to $G = RMult(Q)$ using the subgroup $G_c$ for any $c \in Q$ and coset representatives $L$ and $R$, respectively.  Some of these lead to new examples of Construction \ref{const2}, where $G_c$ does not have a complement, as we now illustrate.
\vspace{1ex}

\textbf{Example 1}  Table \ref{loop} defines a quasigroup $Q_6$ since the table is visibly a bordered Latin square.  We calculate $\lambda_1 = (1)$, $\lambda_2 = (123)(456)$, $\lambda_3 = (132)(465)$, $\lambda_4 = (14)(2536)$, $\lambda_5 = (15)(2634)$, and $\lambda_6 = (16)(2435)$, all even permutations.

\begin{table}[h!]
\begin{center}
\begin{tabular}{c||cccccc}
\null & 1 & 2 & 3 & 4 & 5 & 6 \\ \hline \hline
1 & 1 & 2 & 3 & 4 & 5 & 6 \\
2 & 2 & 3 & 1 & 5 & 6 & 4 \\
3 & 3 & 1 & 2 & 6 & 4 & 5 \\
4 & 4 & 5 & 6 & 1 & 3 & 2 \\
5 & 5 & 6 & 4 & 2 & 1 & 3 \\
6 & 6 & 4 & 5 & 3 & 2 & 1 \\
\end{tabular}
\caption{Quasigroup $Q_6$} \label{loop}
\end{center}
\end{table}

Therefore, $G = LMult(Q_6)$ is a subgroup of the alternating group $A_6$.  Suppose $G_1$ has a complement $C$. Then $C$ is a regular subgroup by Proposition \ref{permlemma}, contradicting Proposition \ref{an}.   Thus, applying Construction 2L (and Proposition \ref{need1}) to G using the subgroup $G_1$ and the left coset representatives $L = \{\lambda_1, \lambda_2, \lambda_3, \lambda_4, \lambda_5, \lambda_6 \}$ gives a new instance of the construction.

By the way, GAP [\ref{gap}] tells us $G$ has order 36 and \\ $G_1 = \{(1), (456), (465), (23)(56), (23)(45), (23)(46) \}$.  The corresponding Cayley-Sudoku Table from Construction 2L has blocks of dimension $6 \times 6$.

On the other hand, $\rho_1 = (1)$, $\rho_2 = (123)(456)$, $\rho_3 = (132)(465)$, $\rho_4 = (14)(25)(36)$, $\rho_5 = (15)(26)(34)$, and $\rho_6 = (16)(24)(35)$ are not all even permutations.  Using GAP again, $G =RMult(Q_6)$ has order 18 and $G_1 = \{(1),(456),(465)\}$.  The corresponding Cayley-Sudoku Theorem from Construction 2R has blocks of dimension $6 \times 3$.  However, $G_1$ does have a complement, namely $\left< \rho_2 \rho_4 \right> = \left< (153426) \right>$, so, unfortunately, this does not give a new instance of Construction 2.

\vspace{1ex}

\textbf{Example 2}  Example 1 readily generalizes to any quasigroup $Q_n$ of order $n>2$ where $n \equiv 2 (\textrm{mod } 4)$ and the left (or right) translations are even permutations.  Table \ref{geneg} illustrates just such a generalization. To verify it is a quasigroup, consider the four subtables formed by the dashed lines.  The upper left and lower right subtables contain the numbers 1 through $\frac{n}{2}$ while the other two contain the numbers $\frac{n}{2}+1$ through $n$.  Moreover, each successive row in the lower right subtable is the previous row shifted (with wrapping) one place to the right.  In the remaining subtables, rows are shifted one place to the left.  Thus we see the table is a bordered Latin square.

One can calculate $\lambda_i$ to be the permutation $\left( (1,2, \ldots, \frac{n}{2})(\frac{n}{2}+1,\frac{n}{2}+2, \ldots, n)\right)^{i-1}$ for $1 \leq i \leq \frac{n}{2}$, and $(1, i)(2, i+1, 3, i+2, \ldots, \frac{n}{2}-1, i+\frac{n}{2}-2, \frac{n}{2}, i+\frac{n}{2}-1)$ for $\frac{n}{2}+1 \leq i \leq n$ where addition is done modulo $n/2$.  It is not hard to see each $\lambda_i$ is an even permutation under the hypotheses on $n$.

\begin{table}[h!]
\begin{tiny}
\begin{center}
\begin{tabular}{c||cccccc:cccccc}
  \null & 1 & 2 & 3 & $\cdots$ & $\frac{n}{2} - 1$ & $\frac{n}{2}$
  & $\frac{n}{2} + 1$ & $\frac{n}{2} + 2$ & $\frac{n}{2} + 3$ & $\cdots$ & $n-1$ & $n$  \\
  \hline \hline
  1 & 1 & 2 & 3 & $\cdots$ & $\frac{n}{2} - 1$ & $\frac{n}{2}$ & $\frac{n}{2} + 1$
  & $\frac{n}{2} + 2$ & $\frac{n}{2} + 3$ & $\cdots$ & $n-1$ & $n$  \\
  2 & 2 & 3 & 4 & $\cdots$ &  $\frac{n}{2}$ & 1  & $\frac{n}{2} + 2$ &
  $\frac{n}{2} + 3$ & $\frac{n}{2} + 4$ & $\cdots$ & $n$ & $\frac{n}{2} + 1$ \\
  3 & 3 & 4 & 5 & $\cdots$ &  1 & 2
  & $\frac{n}{2} + 3$ & $\frac{n}{2} + 4$ & $\frac{n}{2} + 5$ & $\cdots$ & $\frac{n}{2} + 1$ & $\frac{n}{2} + 2$ \\
  $\vdots$ & $\vdots$ & $\vdots$ & $\vdots$ & & $\vdots$ & $\vdots$ & $\vdots$ & $\vdots$ & $\vdots$ & & $\vdots$ & $\vdots$ \\
  $\frac{n}{2}$ & $\frac{n}{2}$ & 1 & 2 & $\cdots$ & $\frac{n}{2} - 2$ & $\frac{n}{2} - 1$
  & $n$ & $\frac{n}{2} + 1$ & $\frac{n}{2} + 2$ & $\cdots$ & $n-2$ & $n-1$ \\
  \hdashline
  $\frac{n}{2} + 1$ & $\frac{n}{2} + 1$ & $\frac{n}{2} + 2$ & $\frac{n}{2} + 3$ & $\cdots$ & $n-1$ & $n$
  & 1 & 3 & 4 & $\cdots$ & $\frac{n}{2}$ & 2 \\
  $\frac{n}{2} + 2$  & $\frac{n}{2} + 2$ & $\frac{n}{2} + 3$ & $\frac{n}{2} + 4$ & $\cdots$ & $n$ & $\frac{n}{2} + 1$
  & 2 & 1 & 3 & $\cdots$ & $\frac{n}{2} - 1$ & $\frac{n}{2}$ \\
  $\frac{n}{2} + 3$ & $\frac{n}{2} + 3$ & $\frac{n}{2} + 4$ & $\frac{n}{2} + 5$ & $\cdots$ & $\frac{n}{2} + 1$ & $\frac{n}{2} + 2$
  & $\frac{n}{2}$ & 2 & 1 & $\cdots$ & $\frac{n}{2} - 2$ & $\frac{n}{2} - 1$ \\
   $\vdots$ & $\vdots$ & $\vdots$ & $\vdots$ & & $\vdots$ & $\vdots$ & $\vdots$ & $\vdots$ & $\vdots$ & & $\vdots$ & $\vdots$ \\
  $n$ & $n$ & $\frac{n}{2} + 1$ & $\frac{n}{2} + 2$ & $\cdots$ & $n-2$ & $n-1$
  & 3 & 4 & 5 & $\cdots$ &  2 & 1 \\
  \hline
\end{tabular}
\caption{Quasigroup $Q_n$} \label{geneg}
\end{center}
\end{tiny}
\end{table}

 As a matter of fact, as long as the left or right translations lie in any permutation group not having a regular subgroup, no complement exists.  It is interesting to know such groups exist of every order except 1, a prime, or a prime squared [\ref{xu}, Theorem 4].

\paragraph*{Recent Appearances of Construction \ref{const1}}
Construction \ref{const1} continues to be rediscovered in special cases.  We find it, for example, in [\ref{pv}] and [\ref{lorch}].  In this section, we briefly sketch the ideas.  A determined reader can fill-in the details.

Let $K$ be a finite field of order $q^2$, $F$ its subfield of order $q$ viewed as a subgroup of $K$ under addition, and $c_0, c_1, \ldots, c_{q-1}$ a complete set of (right)\footnote{The parentheses remind us that the distinction between right and left cosets is unimportant in the abelian group $K$.} coset representatives of $F$ in $K$.  For each $x \in K \backslash F$, it is easy to prove the sets $(F+c_0)x$, \ldots, $(F+c_{q-1})x$ partition $K$ into complete sets of (left) coset representatives of $F$ in $K$.  Although they describe it differently, Pedersen and Vis [\ref{pv}] form the table $L_x$ by arranging the addition table for $K$ as
\begin{center}
\begin{tabular}{c||c|c|c|c|}
\null & $F+c_0$ & $F+c_1$ & \ldots & $F+c_{q-1}$ \\ \hline \hline
$(F+c_0)x$ & & & & \\ \hline
$(F+c_1)x$ & & & & \\ \hline
\vdots & & & & \\ \hline
$(F+c_1)x$ & & & & \\ \hline
\end{tabular}.
\end{center}

We recognize $L_x$ as a Cayley-Sudoku table for $K$ produced using the subgroup $F$ and Construction 1R.

(In pursuit of their main objective, Pedersen and Vis also show $L_x$ and $L_y$ are orthogonal Latin squares when $x \neq y$, thus giving examples of orthogonal Cayley-Sudoku tables. Searching for other sets of mutually orthogonal Cayley-Sudoku tables is interesting.)

Finding Construction 1 in [\ref{lorch}] is more involved.  Again take a finite field $F = \{a_0, a_1, \ldots, a_{q-1}\}$ of order $q$.  Lorch [\ref{lorch}] studies ``sudoku solutions of order $q^2$," that is, $q^2 \times q^2$ Latin squares on the symbols $0,1, \ldots, q^2-1$ subdivided into $q \times q$ blocks with each of the symbols appearing once in each block--just as one would expect when generalizing from sudoku puzzles where $q=3$ (and the symbol $9$ is used instead of $0$).  He uses an addressing scheme to associate each entry in a sudoku solution with an element of $F^4$ and looks, in particular, at sudoku solutions where, for some 2-dimensional subspace $g$ of the vector space $F^4$, the set of addresses for each symbol form a coset of $g$ in $F^4$. (These are called parallel linear sudoku solutions.)  Take any such sudoku solution (e.g. Figure 1 or 4 of [\ref{lorch}]), border it on the left and on the top with a copies of the first column and the first row, respectively, and then replace each symbol in the resulting table with the coset of $g$ containing its addresses.  With some work, one may prove the table so constructed is precisely the Cayley-Sudoku table for the factor group $F^4/g$ under addition obtained from Construction 1L as
\begin{center}
\begin{tabular}{c||c|c|c|c|}
\null & $R_0$ & $R_1$ & \ldots & $R_{q-1}$ \\ \hline \hline
$(a_0,0,0,0)+S$ & & & & \\ \hline
$(a_1,0,0,0)+S$ & & & & \\ \hline
\vdots & & & & \\ \hline
$(a_{q-1},0,0,0)+S$ & & & & \\ \hline
\end{tabular}
\end{center}
where $S:=\{(0,b,0,0)+g : b \in F\}$ is a subgroup of $F^4/g$ and $R_i :=\{(0,0,a_i,d)+g : d \in F\}$ for $0 \leq i \leq q-1$ are complete sets of (right) coset representatives of $S$ in $F^4/g$.  (Lorch also produces what we would call orthogonal Cayley-Sudoku tables and his can be made ``magic'' in the sense of a magic square.)

\paragraph*{The Heritage of Construction 3?}  Construction 3, dubbed the ``centerpiece'' of [\ref{csw}] by a referee, shows how to induce a Cayley-Sudoku table of a group from a Cayley-Sudoku table of any of its subgroups.  So far, we have not encountered other incarnations of Construction 3.  Can any reader find its heritage?

\paragraph*{Acknowledgment}
Portions of this paper come from the first author's honors thesis at Western Oregon University, supervised by the second author.

\subsection*{References}
\begin{enumerate}
\item \label{baer} R. Baer, Nets and Groups, \emph{Trans. Amer. Math. Soc.} \textbf{46} (1939) 110-141.

\item \label{csw}  J. Carmichael, K. Schloeman, and M. B. Ward, Cosets and Cayley-Sudoku Tables, this \textsc{ Magazine}, \textbf{83} (2010) 130-139.

\item \label{denes1}  J. D\'{e}nes, Algebraic and Combinatorial Characterizations of Latin Squares I, \emph{Math. Slovaca} \textbf{17} (1967) 249-265.

\item \label{denes2}  J. D\'{e}nes and A. D. Keedwell, \emph{Latin Squares and Their Applications}, Academic Press, New York and London, 1974.

\item \label{gap} The GAP Group, GAP -- Groups, Algorithms, and Programming, Version 4.6.2 \texttt{http://www.gap-system.org} 2013.

\item \label{lorch}  J. Lorch, Magic Squares and Sudoku, \emph{Amer. Math. Monthly} \textbf{119} (2012) 759-770.

\item \label{pv} R. M. Pedersen and T. L. Vis, Sets of Mutually Orthogonal Sudoku Latin Squares, \emph{College Math. J} \textbf{40} (2009) 174-180.

\item \label{xu} M-Y Xu, A Note on Permutation Groups and Their Regular Subgroups, \emph{J. Aust. Math. Soc.} \textbf{85} (2008) 283-287.

\end{enumerate}

\end{document}